\pdfoutput=1
\RequirePackage{silence}
\documentclass[a4paper,svgnames,11pt]{amsart}
\usepackage[hmarginratio={1:1},vmarginratio={1:1},lmargin=60.0pt,tmargin=60.0pt]{geometry}

\usepackage{booktabs}
\allowdisplaybreaks

\usepackage[T1]{fontenc}
\usepackage{latexsym,exscale,amsfonts,amssymb,mathtools}
\usepackage{amsmath,amsthm,amsfonts,amssymb,amscd,textcomp,bbm}
\usepackage{mathrsfs,stackrel}
\usepackage{etoolbox}
\usepackage{bbm}
\usepackage{mathabx}
\usepackage{latexsym}
\usepackage{graphicx}
\usepackage{subfiles}
\usepackage{xcolor}
\usepackage{xfrac}

\usepackage{tikz}
\usetikzlibrary{decorations}
\usetikzlibrary{decorations.markings}
\usetikzlibrary{decorations.pathreplacing}
\usetikzlibrary{decorations.pathmorphing}
\tikzset{anchorbase/.style={baseline={([yshift=-0.5ex]current bounding box.center)}},
}
\tikzstyle directed=[postaction={decorate,decoration={markings,mark=at position #1 with {\arrow[line width=0.3mm, black]{>}}}}]

\usepackage{enumitem}
\setlist[enumerate]{itemsep=0.15cm,label=\emph{\upshape(\alph*)}}
\setlist[enumerate,2]{itemsep=0.15cm,label=\emph{\upshape(\roman*)}}
\setlist[enumerate,3]{itemsep=0.15cm,label=\emph{\upshape(\Alph*)}}

\let\emph\relax
\DeclareTextFontCommand{\emph}{\bfseries\em}

\newcommand{\N}{{\mathbb{Z}}_{\geq 0}}
\newcommand{\C}{{\mathbb{C}}}

\newcommand{\K}{{\mathbb{K}}}

\newcommand{\R}{{\mathbb{R}}}

\newcommand{\Z}{{\mathbb{Z}}}

\renewcommand{\phi}{\varphi}

\usepackage{xcolor,colortbl}
\definecolor{orchid}{RGB}{143,40,194}
\definecolor{lava}{RGB}{207,16,32}
\definecolor{mydarkblue}{RGB}{10,10,170}
\definecolor{tomato}{RGB}{255,99,71}

\usepackage{aliascnt}

\def\NewTheorem#1{%
\newaliascnt{#1}{equation}%
\newtheorem{#1}[#1]{#1}%
\aliascntresetthe{#1}%
\expandafter\def\csname #1autorefname\endcsname{#1}%
}
\def\equationautorefname~#1\null{(#1)\null}

\numberwithin{equation}{subsection}

\NewTheorem{Proposition}
\NewTheorem{Theorem}
\NewTheorem{Corollary}
\AtEndEnvironment{Corollary}{\null\hfill$\square$}%
\NewTheorem{Lemma}
\theoremstyle{definition}
\NewTheorem{Definition}
\AtEndEnvironment{Definition}{\null\hfill$\Diamond$}%
\NewTheorem{Classification Problem}
\AtEndEnvironment{Classification Problem}{\null\hfill$\Diamond$}%
\NewTheorem{Notation}
\AtEndEnvironment{Notation}{\null\hfill$\Diamond$}%
\NewTheorem{Example}
\AtEndEnvironment{Example}{\null\hfill$\Diamond$}%
\NewTheorem{Examples}
\AtEndEnvironment{Examples}{\vskip-10mm\null\hfill$\Diamond$}%
\NewTheorem{Conjecture}

\theoremstyle{remark}
\NewTheorem{Remark}
\AtEndEnvironment{Remark}{\null\hfill$\Diamond$}%
\NewTheorem{Question}
\AtEndEnvironment{Question}{\null\hfill$\Diamond$}%
\NewTheorem{Assumption}

\usepackage[hypertexnames=false]{hyperref}
\usepackage{bookmark}
\hypersetup{
pdftoolbar=true,
pdfmenubar=true,
pdffitwindow=false,
pdfstartview={FitH},
pdftitle={Growth problems of quantum groups},
pdfauthor={Jensen O'Sullivan and Daniel Tubbenhauer},
pdfsubject={},
pdfcreator={Jensen O'Sullivan and Daniel Tubbenhauer},
pdfproducer={Jensen O'Sullivan and Daniel Tubbenhauer},
pdfkeywords={},
pdfnewwindow=true,
colorlinks=true,
linkcolor=mydarkblue,
citecolor=teal,
filecolor=magenta,
urlcolor=orchid,
linkbordercolor=lava,
citebordercolor=teal,
urlbordercolor=orchid,
linktocpage=true
}

\usepackage{etoolbox}
\def\makeautorefname#1#2{\csdef{#1autorefname}{#2}}
\makeautorefname{section}{Section}%
\makeautorefname{subsection}{Section}%
\makeautorefname{subsubsection}{Section}%

\begin{document}
\title[Growth problems of quantum groups]{Growth problems of quantum groups}
\author[J. O'Sullivan and D. Tubbenhauer]{Jensen O'Sullivan and Daniel Tubbenhauer}

\address{J.O.: The University of Sydney, School of Mathematics and Statistics F07, NSW 2006, Australia}
\email{J.OSullivan@maths.usyd.edu.au}

\address{D.T.: The University of Sydney, School of Mathematics and Statistics F07, Office Carslaw 827, NSW 2006, Australia, \href{http://www.dtubbenhauer.com}{www.dtubbenhauer.com}, \href{https://orcid.org/0000-0001-7265-5047}{ORCID 0000-0001-7265-5047}}
\email{daniel.tubbenhauer@sydney.edu.au}

\begin{abstract}
We study the asymptotic size of decompositions of tensor powers of tilting modules for quantum groups (mostly at a complex root of unity). In type A1 we obtain a sharp result for the number of indecomposable summands, explained by a one dimensional half-line random walk with a periodic congruence constraint. In general type we prove a universal law: the dominant part is governed only by the dimension of the module, while the correction depends only on the root system, so the asymptotic size is largely independent of the specific tilting module.
\end{abstract}


\addtocontents{toc}{\protect\setcounter{tocdepth}{1}}

\maketitle

\tableofcontents

\maketitle

\section{Introduction}\label{S:Intro}

To set the stage, let
\[
b_n=b_n(T)=\#\{\text{indecomposable summands of }T^{\otimes n}\ \text{counted with multiplicity}\},
\]
for an object $T$ in a monoidal category where this count makes sense (e.g. a finite dimensional representation of a group $G$).
Exact multiplicities are extremely sensitive to the precise setting and often completely out of reach, but $b_n$ often has a stable, representation-theoretically natural asymptotic, and even a general form. 
As an analogy: exact prime counts are subtle (at best), yet their growth is governed by the simple formula $n/\ln n$; likewise, our focus is the coarse law for $b_n$.

While writing this paper, the systematic study of asymptotics for the number of indecomposable summands in tensor powers is rather recent, and a general picture has only crystallized in the last few years. 
Following \cite{CoEtOsTu-growth-fractal}, and writing $\sim$ for asymptotic equality, one frequently sees:
\begin{gather}\label{Eq:MainEq}
\scalebox{1.1}{\fcolorbox{black}{orchid!2}{$b_{n}\sim
h(n)\cdot n^{\tau}\cdot\beta^{n},
\quad
\begin{aligned}
&h\colon\N\to\R_{>0}\text{ is a function \emph{bounded away from $0,\infty$}},
\\[-0.1cm]
&n^{\tau}\text{ is the \emph{subexponential factor}, $\tau\in\R$},
\\[-0.1cm]
&\beta^{n}\text{ is the \emph{exponential factor}, 
$\beta\in\R_{\geq 1}$}.
\end{aligned}$}}
\end{gather}
Here bounded away from $0$ and $\infty$ means $c\leq h(n)\leq C$ for $c,C\in\R_{>0}$.

\begin{Remark}
The formula \autoref{Eq:MainEq} is restrictive: it forces exponential growth with a polynomial-type subexponential factor. In particular, terms of the form \(2^{\sqrt{n}}\cdot\beta^n\) (whose subexponential part is not polynomially bounded) or \(n^{(-1)^n}\cdot\beta^n\) (polynomially bounded but not of the form $n^\tau$) cannot occur.

There are, however, natural examples where \autoref{Eq:MainEq} fails completely. For the partition category with its generating object one has the asymptotic
$
b_n \sim
(2n/W(2n))^n\cdot
\exp(n/W(2n) + (2n/W(2n))^{\!1/2} - n - \frac{7}{4})/ \sqrt{1 + W(2n)}
$,
where \(W\) is the Lambert \(W\)-function, and $\sqrt[n]{b_n}\sim\frac{2}{e}\cdot\frac{1}{\log 2n}\cdot n$.

Empirically, \autoref{Eq:MainEq} seems to characterize group-like situations.
\end{Remark}

In practice, the exponential base $\beta$ is $\dim_\K T$ \cite{CoOsTu-growth} (hence independent of $G$, and only partially dependent on $T$), and $h(n)$ is constant or alternates among finitely many constants. Moreover, for reductive groups over $\C$, the subexponential exponent $\tau$ is determined by the root system (thus independent of $T$) \cite{Bi-asymptotic-lie,CoEtOs-growth}. 
This principle does admit exceptions, but it holds in many natural settings, and the surrounding literature has grown rapidly in recent years; see, e.g., the works cited above and related developments \cite{BeSy-non-projective-part,GrTu-diagram-growth,He-group-growth,HeTu-monoid-growth,HeTu-tensor-monoid,KhSiTu-monoidal-cryptography,LaTuVa-growth-pfdim,LaTuVa-growth-pfdim-inf,LaPoReSo-growth-qgroup,Sh-trivial-summands,VaZh-p-groups}.

Let $p$ be the characteristic of the underlying field.
Together these advances support the heuristic
\emph{``exponential from dimension, and subexponential from root system and $p$''}
as an organizing principle for reductive-group-like settings across semisimple and nonsemisimple contexts. We do not know in what generality this works, but:
In characteristic zero, the number of constituents in $T^{\otimes n}$ has exponential base $\beta=\dim_\C T$ with a universal subexponential factor $\tau$ equals minus $1/2$ times the number of positive roots. 
In positive characteristic not much is known currently, but type~A1 results support this heuristic \cite{CoEtOsTu-growth-fractal}.

Our results place the root-of-unity quantum case in this landscape, and show it is close to characteristic~0, and the heuristic works, in particular, there is \emph{no dependence of the quantum characteristic}.

More precisely, the contributions of this paper are as follows.
Let $U_q(\mathfrak{g})$ be a divided power quantum group at a complex root $q$ of unity (e.g. $q=\exp(\pi i/\ell)$) and $T$ a tilting $U_q(\mathfrak{g})$-module. Let $b_n=b_n(T)$ and let $a_n$ be the (quantum) characteristic zero analog.
We study the growth of $b_n$ as $n\to\infty$, and 
show that it is asymptotically of the form $\Theta(a_n)$ (meaning $b_n\in O(a_n)$ and $a_n\in O(b_n)$).
Equivalently, $b_n$ satisfies \autoref{Eq:MainEq} with the same $\tau,\beta$ as in characteristic zero, but we do not
control the finer oscillatory factor $h(n)$ in general.
In type A1 (i.e.\ $\mathfrak{g}=\mathfrak{sl}_2$) and $T$ the defining representation, we go further and determine $h(n)$ explicitly: it is constant for odd $\ell$ and two-periodic for even $\ell$.

Our main results are:
\begin{enumerate}
\item \textbf{Type A1: sharp asymptotics.}

\noindent For the two dimensional indecomposable tilting object, $b_n$ is encoded by a half-line random walk with a periodic (mod~$\ell$) constraint. 
A kernel-method analysis of a bivariate generating function yields a precise local expansion and an explicit leading constant. 
We also comment on the mixed case.
\item \textbf{General simple type: a universal law.}

\noindent For any tilting $T$, the exponential growth rate $\beta$ of $b_n$ equals $\dim_\C T$, and the subexponential exponent $\tau$ depends only on the number of positive roots.
\end{enumerate}
Finally, the random-walk model for type~A1 is interesting from a purely analytic-combinatorics viewpoint (see \autoref{S:SL2}); that section is written so as to require no prior knowledge of representation theory.
\medskip

\noindent\textbf{Acknowledgments.}
We would like to thank the incredible Kevin Coulembier, whose contribution to this work can hardly be overestimated, and we are grateful to Henning Haahr Andersen for suggesting a strategy for the proof of \autoref{L:Bound}.
No other collaborators were harmed in the making of this paper.

The paper is part of the PhD research of the first author, supervised by Kevin Coulembier.
DT subscribes to ``I come from nothing and I go back to nothing. What have I lost?'' and acknowledges support from ARC Future Fellowship FT230100489.

\section{The half-line random walk and a variation}\label{S:SL2}

We will study the following two sequences and relate them to one another. We explain this in a purely analytic combinatorial way and the reader does not need any background in representation theory for this part (knowing standard techniques in analytic combinatorics is enough, as in e.g. \cite{FlSe-analytic-combinatorics}); details on the representation theory are given in \autoref{S:Background}.

We define two sequences $a_n$ and $b_n$, the former being the classical (simple and symmetric) random walk on the half-line (i.e. a hard wall at $-1$), and the latter is a peculiar modular version of it.
Specifically, for $n,m\in\N$, and $\ell\in\Z_{\geq 2}$ let
\begin{gather}\label{Eq:Recursion}
\begin{gathered}
a_{n,m}=\left\{\#\text{ of unit step paths of length $n$ on $\N$ starting at $0$ and ending at $m$},\right.
\\
b_{n,m}=b_{n,m}(\ell)=
\begin{cases}
a_{n,m} & \text{ if } m \equiv -1\bmod\ell,\\
b_{n-1,m-1} + b_{n-1,m+1} & \text{ if } m \equiv m_0\bmod\ell \text{ where } 0 \leq m_0 < \ell -2,\\
b_{n-1,m-1} & \text{ if } m \equiv - 2 \bmod\ell.
\end{cases}
\end{gathered}
\end{gather}
Note that, if $\ell=2$, then there is no middle case for $b_{n,m}$.
We will study the following two sequences:
\begin{gather*}
a_n=\sum_{m\geq 0} a_{n,m},\quad
b_n=\sum_{m\geq 0} b_{n,m}.
\end{gather*}
(Note that, in this case as well as below, summing over $m\geq 0$ is the same as summing over $0\leq m\leq n$).
The study of these is motivated by \autoref{S:SL2one}, but we do not need to know this for now.

\begin{Remark}
Strictly speaking $a_n$ and $b_n$ are scaled versions of random walks.
\end{Remark}

\begin{Example}
We can identify $a_{n}$ with words of length $n$ in the alphabet $\{+1,-1\}$ (or just $\{+,-\}$) such that partial sums are all nonnegative, and $a_{n,m}$ are the words with final sum $m$. For example,
\begin{gather*}
++-+
\leftrightsquigarrow
\begin{tikzpicture}[anchorbase,scale=1]
\draw[ultra thick] (-1.5,0) to (2.5,0);
\node at (-2,0){$\dots$};
\draw[thick,fill=tomato] (-1,0) circle (0.1cm)node[below,yshift=-0.05cm]{$-1$};
\draw[thick,fill=tomato] (0,0) circle (0.1cm)node[below,yshift=-0.05cm]{$\phantom{1}0\phantom{1}$};
\draw[thick,fill=tomato] (1,0) circle (0.1cm)node[below,yshift=-0.05cm]{$1$};
\draw[thick,fill=tomato] (2,0) circle (0.1cm)node[below,yshift=-0.05cm]{$2$};
\node at (3,0){$\dots$};
\draw[directed=0.99] (0,0.25) to (1,0.25);
\draw[directed=0.99] (1,0.25) to (2,0.25);
\draw[directed=0.99] (2,0.5) to (1,0.5);
\draw[directed=0.99] (1,0.75) to (2,0.75);
\end{tikzpicture},
\end{gather*}
is one of the two paths giving $a_{4,2}=2$.
Let $\ell=4$. Then one then easily counts that
\begin{gather*}
\scalebox{0.5}{$\left (
\begin {array} {cccccccccccccccc} 1 & 0 & 0 & 0 & 0 & 0 & 0 & 0 & 0 \
& 0 & 0 & 0 & 0 & 0 & 0 & 0 \\ 0 & 1 & 0 & 0 & 0 & 0 & 0 & 0 & 0 & 0 \
& 0 & 0 & 0 & 0 & 0 & 0 \\ 1 & 0 & 1 & 0 & 0 & 0 & 0 & 0 & 0 & 0 & 0 \
& 0 & 0 & 0 & 0 & 0 \\ 0 & 2 & 0 & 1 & 0 & 0 & 0 & 0 & 0 & 0 & 0 & 0 \
& 0 & 0 & 0 & 0 \\ 2 & 0 & 3 & 0 & 1 & 0 & 0 & 0 & 0 & 0 & 0 & 0 & 0 \
& 0 & 0 & 0 \\ 0 & 5 & 0 & 4 & 0 & 1 & 0 & 0 & 0 & 0 & 0 & 0 & 0 & 0 \
& 0 & 0 \\ 5 & 0 & 9 & 0 & 5 & 0 & 1 & 0 & 0 & 0 & 0 & 0 & 0 & 0 & 0 \
& 0 \\ 0 & 14 & 0 & 14 & 0 & 6 & 0 & 1 & 0 & 0 & 0 & 0 & 0 & 0 & 0 \
& 0 \\ 14 & 0 & 28 & 0 & 20 & 0 & 7 & 0 & 1 & 0 & 0 & 0 & 0 & 0 & 0 \
& 0 \\ 0 & 42 & 0 & 48 & 0 & 27 & 0 & 8 & 0 & 1 & 0 & 0 & 0 & 0 & 0 \
& 0 \\ 42 & 0 & 90 & 0 & 75 & 0 & 35 & 0 & 9 & 0 & 1 & 0 & 0 & 0 & 0 \
& 0 \\ 0 & 132 & 0 & 165 & 0 & 110 & 0 & 44 & 0 & 10 & 0 & 1 & 0 & 0 \
& 0 & 0 \\ 132 & 0 & 297 & 0 & 275 & 0 & 154 & 0 & 54 & 0 & 11 & 0 \
& 1 & 0 & 0 & 0 \\ 0 & 429 & 0 & 572 & 0 & 429 & 0 & 208 & 0 & 65 & 0 \
& 12 & 0 & 1 & 0 & 0 \\ 429 & 0 & 1001 & 0 & 1001 & 0 & 637 & 0 & 273 \
& 0 & 77 & 0 & 13 & 0 & 1 & 0 \\ 0 & 1430 & 0 & 2002 & 0 & 1638 & 0 \
& 910 & 0 & 350 & 0 & 90 & 0 & 14 & 0 & 1 \\\end {array} \right)$}
,
\scalebox{0.5}{$\left (
\begin {array} {cccccccccccccccc} 1 & 0 & 0 & 0 & 0 & 0 & 0 & 0 & 0 & 0 & 0 & 0 & 0 & 0 & 0 & 0 \\ 0 \
& 1 & 0 & 0 & 0 & 0 & 0 & 0 & 0 & 0 & 0 & 0 & 0 & 0 & 0 & 0 \\ 1 & 0 \
& 1 & 0 & 0 & 0 & 0 & 0 & 0 & 0 & 0 & 0 & 0 & 0 & 0 & 0 \\ 0 & 1 & 0 \
& 1 & 0 & 0 & 0 & 0 & 0 & 0 & 0 & 0 & 0 & 0 & 0 & 0 \\ 1 & 0 & 3 & 0 \
& 1 & 0 & 0 & 0 & 0 & 0 & 0 & 0 & 0 & 0 & 0 & 0 \\ 0 & 1 & 0 & 4 & 0 \
& 1 & 0 & 0 & 0 & 0 & 0 & 0 & 0 & 0 & 0 & 0 \\ 1 & 0 & 9 & 0 & 4 & 0 \
& 1 & 0 & 0 & 0 & 0 & 0 & 0 & 0 & 0 & 0 \\ 0 & 1 & 0 & 13 & 0 & 6 & 0 \
& 1 & 0 & 0 & 0 & 0 & 0 & 0 & 0 & 0 \\ 1 & 0 & 28 & 0 & 13 & 0 & 7 \
& 0 & 1 & 0 & 0 & 0 & 0 & 0 & 0 & 0 \\ 0 & 1 & 0 & 41 & 0 & 27 & 0 \
& 7 & 0 & 1 & 0 & 0 & 0 & 0 & 0 & 0 \\ 1 & 0 & 90 & 0 & 41 & 0 & 34 \
& 0 & 9 & 0 & 1 & 0 & 0 & 0 & 0 & 0 \\ 0 & 1 & 0 & 131 & 0 & 110 & 0 \
& 34 & 0 & 10 & 0 & 1 & 0 & 0 & 0 & 0 \\ 1 & 0 & 297 & 0 & 131 & 0 \
& 144 & 0 & 54 & 0 & 10 & 0 & 1 & 0 & 0 & 0 \\ 0 & 1 & 0 & 428 & 0 \
& 429 & 0 & 144 & 0 & 64 & 0 & 12 & 0 & 1 & 0 & 0 \\ 1 & 0 & 1001 & 0 \
& 428 & 0 & 573 & 0 & 273 & 0 & 64 & 0 & 13 & 0 & 1 & 0 \\ 0 & 1 & 0 \
& 1429 & 0 & 1638 & 0 & 573 & 0 & 337 & 0 & 90 & 0 & 13 & 0 & 1 \\\end {array} \right)$}
,
\end{gather*}
are the tables (rows are indexed by $n$) of the numbers $a_{n,m}$ and $b_{n,m}$, respectively.
\end{Example}

\begin{Example}
For $\ell\in\{2,3,4,5\}$, the first eleven values of $b_n$ are
\begin{gather*}
\ell=2\colon\{1, 1, 1, 3, 3, 10, 10, 35, 35, 126, 126\},\quad
\ell=3\colon\{1, 1, 2, 2, 5, 6, 15, 21, 50, 77, 176\},\\
\ell=4\colon\{1, 1, 2, 3, 5, 9, 14, 29, 43, 99, 142\},\quad
\ell=5\colon\{1, 1, 2, 3, 6, 9, 19, 28, 62, 91, 208\}.
\end{gather*}
These are easily generated from the definition and \autoref{L:BinomialA}. For $n>0$, dividing these by the asymptotic approximation $a_n^a=\sqrt{2/\pi}\cdot n^{-1/2}\cdot 2^n$ of $a_n$ (call the result $b_n^{\prime}=b_n/a_n^a$) and plotting them for $n\in\{1,\dots,100\}$ gives
\begin{gather*}
\begin{tikzpicture}[anchorbase]
\node at (0,0) {\includegraphics[height=5.5cm]{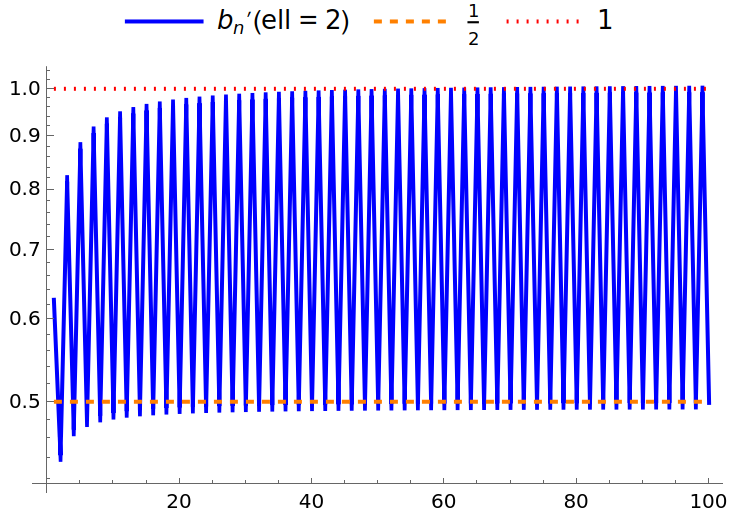}};
\end{tikzpicture},
\begin{tikzpicture}[anchorbase]
\node at (0,0) {\includegraphics[height=5.5cm]{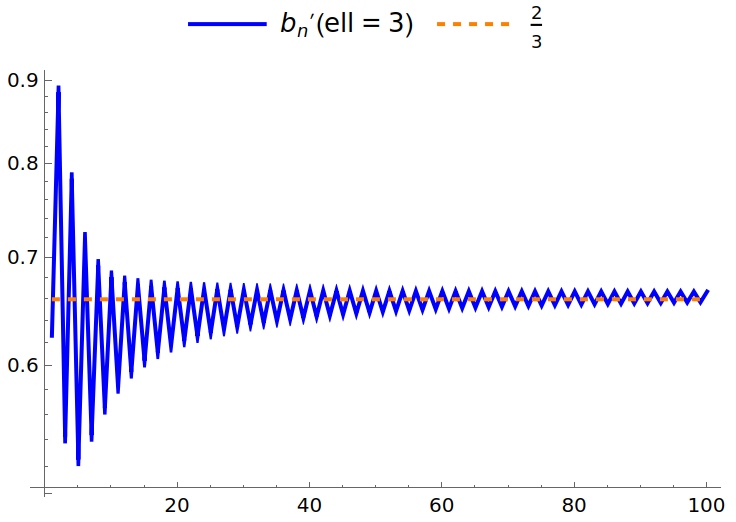}};
\end{tikzpicture},
\\
\begin{tikzpicture}[anchorbase]
\node at (0,0) {\includegraphics[height=5.5cm]{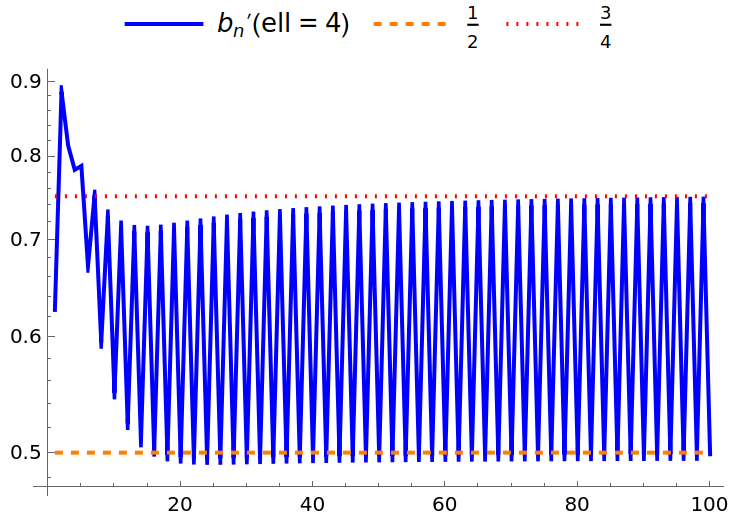}};
\end{tikzpicture},
\begin{tikzpicture}[anchorbase]
\node at (0,0) {\includegraphics[height=5.5cm]{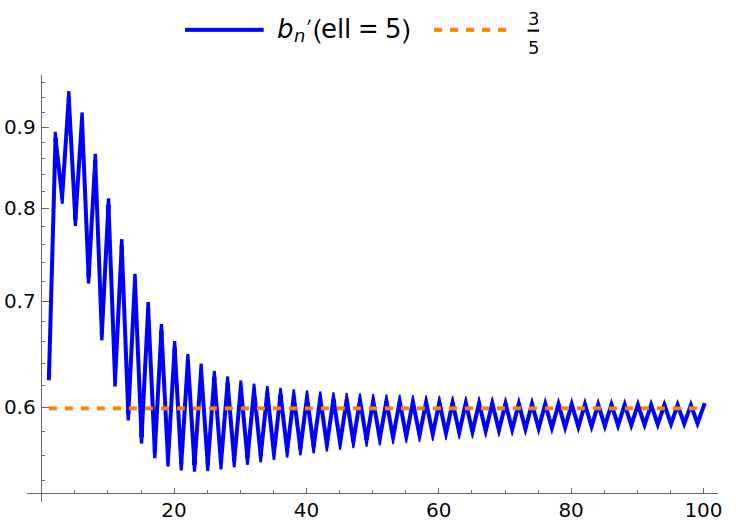}};
\end{tikzpicture}.
\end{gather*}
This motivates \autoref{T:MainSL2} below.
\end{Example}

\begin{Lemma}\label{L:BinomialA}
We have $a_{n,m}=0$ for $m\not\equiv n\bmod 2$, and otherwise $a_{n,m}=\binom{n}{(n-m)/2}-\binom{n}{(n-m)/2-1}$.
\end{Lemma}

\begin{proof}
Well-known and sometimes called the ballot theorem.
\end{proof}

\begin{Example}\label{E:elltwo}
The case $\ell=2$ is somewhat special, so we get it out of the way now. In this case $b_{n,m}=a_{n,m}$ for odd $m$, and $b_{n,m}=a_{n-1,m-1}$ for even $m$. Since $a_{n,m}=0$ for $m\not\equiv n\bmod 2$, we get $b_{n}=a_{n}$ for odd $n$ and $b_{n}=a_{n-1}$ for even $n$.
In other words, depending on parity, $b_n$ is the classical random walk on $\N$, or it is this walk where one omits the final step.
\end{Example}

The next lemma makes the below more explicit. 
Here and below, identities of generating functions are to be interpreted in formal power series; when the one side is analytic (e.g. in a certain domain), its Taylor series agrees term-wise with the other side.

\begin{Lemma}\label{L:NotNeeded}
The generating function $a(x,y)=\sum_{n,m\ge 0}a_{n,m}x^n y^m$ for $a_{n,m}$ is
\begin{gather*}
a(x,y)=\tfrac{1-\tfrac{1-\sqrt{\,1-4x^2\,}}{2xy}}{1-x\,(y+y^{-1})}.
\end{gather*}
\end{Lemma}

\begin{proof}
A step-by-step decomposition on the half-line walk gives the standard kernel equation
\[
\bigl(1-x(y+y^{-1})\bigr)\,a(x,y)=1-x\,y^{-1}a(x,0).
\]
Let $y_\pm(x)=\frac{1\pm\sqrt{1-4x^2}}{2x}$ be the roots of $1-x(y+y^{-1})=0$. Among the two roots,
only \(y_-(x)\in\C[[x]]\) whereas \(y_+(x)\notin\C[[x]]\) since $\sqrt{1-4x^2}$ expands as $1-2x-2x^2-4x^4-\dots$. By usual singularity analysis, we substitue in the smallest root \(y=y_-(x)\) into the kernel equation in \(\C((x))\) annihilates
the left–hand side and yields \(0=1-xy_-^{-1}(x)\,a(x,0)\). Hence $a(x,0)=y_-/x=(1-\sqrt{1-4x^2})/(2x^2)$. Substituting back gives
\[
a(x,y)=\tfrac{1-\tfrac{1-\sqrt{1-4x^2}}{2xy}}{1-x(y+y^{-1})},
\]
as desired.
\end{proof}

We will use the usual big O notation, and $\sim$ means asymptotically equal. We use a superscript $a$ for any asymptotic formula.
For fixed $\ell\ge3$ and a residue $r\in\{0,\dots,\ell-1\}$, set
\[
c_n^{(r)}=\sum_{\substack{m\ge0\\ m\equiv r\bmod \ell}} a_{n,m}.
\]
So $c_n^{(r)}$ counts the number of unit step paths of length $n$ on $\N$ starting at $0$ and ending at a vertex congruent to $r$ modulo $\ell$. The next theorem says that such paths are evenly spread when varying $r$. To state it, recall $g\in O(f)$ denotes ($g\leq C\cdot f$ asymptotically for $C\in\R_{>0}$).

\begin{Theorem}\label{T:AlmostMain}
\leavevmode
\begin{enumerate}
\item We have
\begin{gather*}
a_{n}\sim a_{n}^{a}:=\sqrt{2/\pi}\cdot n^{-1/2}\cdot 2^n
\end{gather*}
and 
$|a_{n}-a_{n}^{a}|/(n^{-1/2}\cdot 2^n)\in O(\tfrac1n)$.

\item For $0\leq r<\ell$ we have
\[
c_n^{(r)} \;=\;
\begin{cases}
\tfrac{1}{\ell}\cdot a_n\bigl(1+O(\tfrac1n)\bigr),
& \ell \text{ odd},\\[2ex]
0, & \ell \text{ even and } r\not\equiv n\pmod 2,\\[1.2ex]
\tfrac{2}{\ell}\cdot a_n\,\bigl(1+O(\tfrac1n)\bigr),
& \ell \text{ even and } r\equiv n\pmod 2.
\end{cases}
\]
\end{enumerate}
\end{Theorem}

\begin{proof}
Let $P_{0,m}^{(n)}$ denote the probability of ending at $m$ after $n$ steps, starting at zero, always staying nonnegative, and with probability $p$ to move right at each step (and probability $1-p$ to move left). To model our problem, we fix $p=1/2$, and see that $a_{n,m} = 2^n P_{0,m}^{(n)}$.
For $u\in\mathbb{C}$, write
\[
A_n^{prob}(u)=\sum_{m\ge0} P_{0,m}^{(n)}\,u^m
\]
for the endpoint generating function, which converges for $|u|<1$ by the usual theory of random walks. Denote by $A_n(u) := 2^n A_n^{prob}(u)$ the scaled version of the generating function which has the number of walks as coefficients, rather than probabilities.\\

\textit{Step 1: Spectral integral for $A_n(u)$.}
Let $U_j(x)$ be the $j$th Chebyshev polynomial of the second kind with $U_{0}(x)=1,U_1(x)=2x$ and 
$U_j(x)=2xU_{j-1}(x)-U_{j-2}(x)$ for $j\in\Z_{\geq 2}$.
For the half-line random walk, \cite{KaMcGr-random-walks} gives
\[
P_{0,m}^{(n)}=\int_{-1}^1 x^n\,U_{m}(x)\,w(x)\,dx,\quad
w(x)=\frac{2}{\pi}\sqrt{1-x^2}.
\]
This is the probability of walking from zero to $m$ in $n$ steps. Multiplying by $u^{m}$ and summing over $m\ge 0$, then using
\[
\sum_{m\ge0} U_m(x)u^m=\frac{1}{1-2ux+u^2},
\]
we obtain
\begin{gather*}
A_n(u)=2^nA_n^p(u)=2^n\int_{-1}^1 \frac{x^n\,w(x)}{1-2ux+u^2}\,dx.
\end{gather*}
With the substitution $x=\cos\theta$ (so $dx=-\sin\theta\,d\theta$ and $w(x)\,dx=\frac{2}{\pi}\sin^2\theta\,d\theta$) this becomes
\begin{gather*}
A_n(u)=\frac{2}{\pi}\int_{0}^{\pi} (2\cos\theta)^n\,
\frac{\sin^2\theta}{\,1-2u\cos\theta+u^2\,}\,d\theta
=
\frac{2}{\pi}\int_{0}^{\pi} (2\cos\theta)^n\,
\frac{\sin^2\theta}{(e^{i\theta}-u)(e^{-i\theta}-u)}\,d\theta,
\end{gather*}
and boundary values on $|u|=1$ are taken by radial limits $u=\rho e^{i\phi}$ with $\rho\uparrow1$.\\

\textit{Step 2: Special values $u=\pm1$.}
At $u=1$ the integrand has a removable singularity:
\[
\frac{\sin^2\theta}{1-2\rho\cos\theta+\rho^2}\xrightarrow[\rho\uparrow1]{}\frac{\sin^2\theta}{2(1-\cos\theta)}
=\cos^2\!\bigl(\frac{\theta}{2}\bigr)=\frac{1+\cos\theta}{2}.
\]
Hence, we get
\begin{gather}\label{eq:An-1}
A_n(1)=\frac{2}{\pi}\int_{0}^{\pi} (2\cos\theta)^n\,\frac{1+\cos\theta}{2}\,d\theta,
\end{gather}
which in turn is equal to $a_n$ by definition.

Since $a_{n,m}=0$ unless $m\equiv n\bmod 2$,
\begin{gather*}
A_n(-1)=\sum_{m\ge0}a_{n,m}(-1)^m=(-1)^n\sum_{m\ge0}a_{n,m}=(-1)^n a_n.
\end{gather*}\\

\textit{Step 3: Asymptotic for $A_n(1)$ and $A_n(-1)$.}
We apply Laplace's method at the endpoint $\theta=0$ (equivalently, Watson's
lemma). Set $\theta=t/\sqrt{n}$. Using
\[
\log(2\cos\theta)=\log 2-\tfrac{\theta^{2}}{2}-\tfrac{\theta^{4}}{12}+O(\theta^{6}),
\quad
\frac{1+\cos\theta}{2}=1-\tfrac{\theta^{2}}{4}+O(\theta^{4}),
\]
one gets
\[
(2\cos\theta)^n
=2^{n}\exp\!\Big(-\tfrac{t^{2}}{2}-\tfrac{t^{4}}{12n}+O(n^{-2})\Big),
\quad
\frac{1+\cos\theta}{2}
=1-\tfrac{t^{2}}{4n}+O(n^{-2}).
\]
Then extend the integral to $\infty$ (the cost is $e^{-cn}$ for some $c>0$,
so we can extend the upper limit to $\infty$ without harm), with $d\theta=dt/\sqrt{n}$ we get:
\[
A_n(1)=\frac{2}{\pi}\frac{2^n}{\sqrt{n}}
\int_0^\infty e^{-t^2/2}\Big(1-\tfrac{t^2}{4n}+O(n^{-2})\Big)\,dt.
\]
Using $\int_0^\infty e^{-t^2/2}dt=\sqrt{\pi/2}$ and
$\int_0^\infty t^2 e^{-t^2/2}dt=\sqrt{\pi/2}$ gives
\begin{gather}\label{eq:An1-asymp}
A_n(1)=\sqrt{2/\pi}\cdot n^{-1/2}\cdot 2^n\bigl(1+O(\tfrac1n)\bigr).
\end{gather}
The case $u=-1$ is analogous, up to an alternating sign.

Thus, the first part of the theorem is proven. It remains to address the second part.\\

\textit{Step 4: Uniform bound for $A_n(u)$ at other roots of unity.}
Fix $\ell\ge3$ and $\omega=e^{2\pi i/\ell}$. For $j\in\{1,\dots,\ell-1\}$ with $\omega^j\neq\pm1$, split $A_n(u)$ at $\delta=n^{-1/3}$:
\[
A_n(\omega^j)=\frac{2}{\pi}\!\left(\int_{0}^{\delta}+\int_{\delta}^{\pi-\delta}
+\int_{\pi-\delta}^{\pi}\right)
(2\cos\theta)^n\,\frac{\sin^2\theta}{\,(e^{i\theta}-\omega^j)(e^{-i\theta}-\omega^j)\,}\,d\theta.
\]
On $[0,\delta]$ we use $\sin\theta\sim\theta$ and $(2\cos\theta)^n\le 2^n e^{-n\theta^2/2}$, while for $\theta$ sufficiently small,
\[
|e^{i\theta}-\omega^j|\ \ge\ |1-\omega^j|-|e^{i\theta}-1|\ \ge\ \tfrac12\,|1-\omega^j|,
\]
and similarly $|e^{-i\theta}-\omega^j|\ge \tfrac12\,|1-\omega^j|$. Hence
\[
\big|(e^{i\theta}-\omega^j)(e^{-i\theta}-\omega^j)\big|\ \ge\ \tfrac14\,|1-\omega^j|^2\ =\ \Delta_j>0,
\]
so
\[
\int_{0}^{\delta}\cdots\ \ll\ \frac{2^n}{\Delta_j}\int_{0}^{\infty}\theta^2 e^{-\frac{n}{2}\theta^2}\,d\theta
= C\cdot n^{-3/2}\cdot 2^n.
\]
The same bound holds on $[\pi-\delta,\pi]$. Put $\vartheta=\pi-\theta$. Then $\sin\theta\approx\vartheta$ and, by the triangle inequality,
\[
|e^{i\theta}+\omega^j|\ \ge\ |{-}1-\omega^j|-|e^{i\theta}+1|\ \ge\ \tfrac12\,|{-}1-\omega^j|
\]
for $\theta$ close to $\pi$, hence
\[
|(e^{i\theta}-\omega^j)(e^{-i\theta}-\omega^j)|\ \ge\ \tfrac14\,|{-}1-\omega^j|^2\ >\ 0.
\]
Therefore the denominator is uniformly bounded away from $0$ while $\sin^2\theta\approx\vartheta^2$, so the integrand is $O(\vartheta^2)$ and the same Laplace–type bound follows.
On $[\delta,\pi-\delta]$ we have $|\cos\theta|\le 1-c\delta^2$, so
$(2\cos\theta)^n\le 2^n e^{-c n\delta^2}=2^n e^{-c n^{1/3}}$, negligible versus $n^{-3/2}\cdot 2^n$.
Therefore, uniformly for such $j$,
\begin{gather}\label{eq:wj-bound-prod}
A_n(\omega^j)=O\!\bigl(n^{-3/2}\cdot 2^n\bigr).
\end{gather}
(When $\ell$ is even and $\omega^{\ell/2}=-1$, we use \eqref{eq:An-1} instead of \eqref{eq:wj-bound-prod}.)\\

\textit{Step 5: Root–of–unity filter.}
For $r\in\{0,\dots,\ell\}$, the root-of-unity filter gives
\[
c_n^{(r)}=\sum_{\substack{m\ge0\\ m\equiv r\bmod \ell}}a_{n,m}
=\frac{1}{\ell}\sum_{j=0}^{\ell-1}\omega^{-rj}\,A_n(\omega^j).
\]
If $\ell$ is odd, then only $j=0$ contributes a main term; by \eqref{eq:An1-asymp} and \eqref{eq:wj-bound-prod},
\[
c_n^{(r)}=\frac{A_n(1)}{\ell}+O\!\bigl(n^{-3/2}\cdot 2^n\bigr)
=\frac{a_n}{\ell}\bigl(1+O(\tfrac1n)\bigr).
\]
If $\ell$ is even, then the indices $j=0$ and $j=\ell/2$ give
\[
\frac{1}{\ell}\bigl(A_n(1)+(-1)^r A_n(-1)\bigr)
=\frac{a_n}{\ell}\bigl(1+(-1)^{r+n}\bigr),
\]
so $c_n^{(r)}=0$ when $r\not\equiv n\bmod 2$ and
$c_n^{(r)}=\frac{2a_n}{\ell}$ when $r\equiv n\bmod 2$, up to the same
$O(n^{-3/2}\cdot 2^n)$ remainder from \eqref{eq:wj-bound-prod}. Normalizing by \eqref{eq:An1-asymp} yields the stated $O(1/n)$ relative error.
\end{proof}

Next, we analyze $b_{n,m}$.
To do this, we fix $\ell$, and find the generating function of the sequence $b_n$, given by 
\[f(x) = \sum_{n \geq 0 } b_n x^n = \sum_{n\geq 0}\sum_{m \geq 0 } b_{n,m} x^n. \]
The recursion \autoref{Eq:Recursion} defining $b_{n,m}$ inspires us to write 
\begin{gather*}
f(x) = \sum_{k=0}^{\ell-1} \sum_{\substack{m \equiv k\bmod\ell\\n \geq k}} b_{n,m}x^n,
\end{gather*}
allowing us to look at the behavior of each $\sum_{m\equiv k\bmod\ell} b_{n,m} x^n$ and sum up multiplicities via congruence class. We do so now.
All congruences below are modulo $\ell$ unless stated otherwise.

\begin{Lemma}\label{L:R Recursions} 
Let $2\leq j \leq \ell$. By setting $R_{\ell - 2}(x) = 1$, and $R_{\ell - 3}(x) = x^{-1}$,
and $R_{\ell-j}(x) = x^{-1} R_{\ell-j+1 }(x) - R_{\ell-j+2}(x)$ for $j>3$, we have
\begin{gather*}
\sum_{\substack{m \equiv \ell - j \\ n\geq \ell-j}} b_{n,m}x^n = R_{\ell - j}(x) \sum_{\substack{m \equiv \ell - 2 \\ n\geq \ell - 2}} b_{n,m}x^n.
\end{gather*}
\end{Lemma}

\begin{proof}
We first verify this for $j\in\{2,3\}$. The case $R_{\ell-2}$ is immediate, and 
\[\sum_{\substack{m \equiv \ell - 2\\n\geq \ell-2}} b_{n,m}x^n = \sum_{\substack{m \equiv \ell - 2\\n\geq \ell-2}} b_{n-1,m-1}x^n\\
= x \sum_{\substack{m \equiv \ell - 3\\n\geq \ell-3}} b_{n,m}x^n.\]
Then, via induction, suppose the result is true for all $k\leq j$, and that $j \leq \ell-1$. Then,
\begin{align*} 
\sum_{\substack{m \equiv \ell - j\\n\geq \ell-j}} b_{n,m}x^n = \sum_{\substack{m \equiv \ell - j\\n\geq \ell-j}} b_{n-1,m-1}x^n + \sum_{\substack{m \equiv \ell - j\\n\geq \ell-j}} b_{n-1,m+1}x^n\\
= x \sum_{\substack{m \equiv \ell - j - 1\\n\geq \ell-j - 1}} b_{n,m}x^n + x\sum_{\substack{m \equiv \ell - j + 1\\n\geq \ell-j + 1}} b_{n,m}x^n\\
= x \sum_{\substack{m \equiv \ell - j - 1\\n\geq \ell-j - 1}} b_{n,m}x^n + x R_{\ell - j + 1}(x) \sum_{\substack{m \equiv \ell - 2 \\ n\geq \ell-2}} b_{n,m} x^n.
\end{align*}
By assumption, \[\sum_{\substack{m \equiv \ell - j \\ n\geq \ell-j}}b_{n,m}x^n = R_{\ell - j}(x) \sum_{\substack{m \equiv \ell - 2 \\ n\geq \ell-2}}b_{n,m}x^n,\] and hence
\begin{gather*}
R_{\ell - j - 1}(x) = \tfrac{R_{\ell - j}(x) - xR_{\ell - j +1}(x)}{x} = x^{-1}R_{\ell - j}(x) - R_{\ell - j +1}(x).
\end{gather*}
We get the stated expression for $j+1$, and the proof is complete.
\end{proof}

\begin{Lemma}\label{L:ChPoly}
For 
$2 \leq j \leq \ell$, we have $R_{\ell-j}(x) = U_{j-2}(1/2x)$. Thus, \[\sum_{\substack{m \equiv \ell-j \\n\geq \ell-j}} b_{n,m}x^n = U_{j-2}(1/2x) \sum_{\substack{m \equiv \ell-2 \\n \geq \ell-2}} b_{n,m}x^n.\]
\end{Lemma}

\begin{proof}
Firstly, $R_{\ell-2}(x) = 1 = U_0(1/2x)$, and $R_{\ell-3}(x) = 1/x = U_1(1/2x)$. Then, by $U_j(x) = 2x U_{j-1}(x) - U_{j-2}(x)$, we have $U_j(1/2x) = x^{-1} U_{j-1}(1/2x) - U_{j-2}(1/2x)$, which is the same recursion as our rational functions satisfy.
\end{proof}

Now, we may write this sum over the first $\ell-2$ congruence classes in terms of the $(\ell-1)$th class. The purpose of this is to then use \autoref{T:AlmostMain} to give us a related asymptotic.

\begin{Lemma}\label{L:RecursionInTermsOfA}
We have that \begin{gather*}
\sum_{\substack{m \equiv \ell - 2\\n\geq \ell-2}} b_{n,m}x^n = \frac{1 + x\sum_{\substack{m \equiv \ell-1\\n \geq \ell -1}} a_{n,m} x^n}{xU_{\ell-1}(1/2x)}.
\end{gather*}
\end{Lemma}

\begin{proof}
The calculation here is much the same as in \autoref{L:R Recursions}. We only need to compare coefficients when splitting up $\sum_{\substack{m \equiv 0\\n \geq 0}} b_{n,m}x^n$ using the rule in (2.1).
\end{proof}

Now, we may use the above results to rewrite our generating function. 

\begin{Lemma}\label{P:GenForB}
The generating function for the sequence $b_n$ is 
\[ \frac{\sum_{k=0}^{\ell-1} U_k(\frac{1}{2x})}{U_{\ell-1}(\frac{1}{2x})}\sum_{\substack{m\equiv \ell-1\\n\geq \ell-1}}a_{n,m}x^n + \frac{ \sum_{k=0}^{\ell-2} U_k(\frac{1}{2x})}{xU_{\ell-1}(\frac{1}{2x})}.\]
\end{Lemma}

\begin{proof}
We use \autoref{L:ChPoly}
and \autoref{L:RecursionInTermsOfA}.
\end{proof}

The following is main theorem of this section.

\begin{Theorem}\label{T:MainSL2}
We have
\[
b_n\sim
b_n^a:=
\begin{cases}
\tfrac{\ell+1}{2\ell} \cdot a_n^a
& \ell \text{ odd},\\[2ex]
\tfrac{\ell+1+(-1)^{n+1}}{2\ell}\cdot a_n^a, & \ell \text{ even, }
\end{cases}
\] 
and 
$|b_{n}-b_{n}^{a}|/(n^{-1/2}\cdot 2^n)\in O(\tfrac1n)$.
\end{Theorem}  

\begin{proof}
We prove this for $\ell$ even. If $\ell$ is odd, it is easier, and almost identical, except there is no $A_n(-1)$ sum to worry about.

For $\ell=2$, \autoref{E:elltwo} and \autoref{T:AlmostMain} give
\begin{gather*}
b_n\sim \tfrac{3+(-1)^{n+1}}{4}\cdot a_n\sim 
\tfrac{3+(-1)^{n+1}}{4}\cdot\sqrt{2/\pi}\cdot n^{-1/2}\cdot 2^n.
\end{gather*}
So let us now assume $\ell\neq 2$.

Recall the generating function for $b_n$ as in \autoref{P:GenForB}.
We may rewrite our generating function using the root-of-unity filter for $c_n^{(\ell-1)}$ from the proof of \autoref{T:AlmostMain} as
\begin{gather*} f(x) = \frac{1}{\ell} \cdot \frac{\sum_{k=0}^{\ell-1} U_k(1/2x)}{U_{\ell-1}(1/2x)} \sum_{n\geq 0} a_n x^n + \frac{1}{\ell} \cdot\frac{\sum_{k=0}^{\ell-1} U_k(1/2x)}{U_{\ell-1}(1/2x)} \sum_{n\geq 0} (-1)^{n+1} a_n x^n \\ + \frac{1}{\ell} \cdot \frac{\sum_{k=0}^{\ell-1} U_k(1/2x)}{U_{\ell-1}(1/2x)} \sum_{\substack{j=1\\j\neq \ell/2\\n\ge 0}}^{\ell-1} \omega^{j(1-\ell)}A_n(\omega^j)x^n +\frac{\sum_{k=0}^{\ell-2} U_k(1/2x)}{xU_{\ell-1}(1/2x)}.
\end{gather*}
Standard singularity analysis tells us we only need to look at the smallest singularities of this function to determine the growth of our sequence $b_n$. 

We already know from \autoref{T:AlmostMain} that the growth of the coefficients of the generating function under consideration is dominated by the first two terms (evaluation at $1$ and $-1$). To easily see the singularities of $\sum_{n\geq 0}a_n x^n$, we relate to the classical random walk case. Recall that $a_n = 2^n \sum_{m\geq 0} P_{0,m}^{(n)}.$ Then, we have 
\[\sum_{n\geq 0} a_n x^n = \sum_{n\ge 0} P_{0,m}^{(n)} (2x)^n, \]
where the singularity in the right-hand side is known to be at $2x =1$ (cf. \autoref{L:NotNeeded}). Similarly, 
\[\sum_{n\geq 0}(-1)^{n+1}a_n x^n = -\sum_{x\geq 0} a_n (-x)^n = -\sum_{x\geq 0} P_n (-2x)^n, \]
giving us a singularity in the second term at $x=-1/2$.

Being a rational function, the potential singularities of \[\frac{\sum_{k=0}^{\ell-1} U_k(\frac{1}{2x})}{U_{\ell-1}(\frac{1}{2x})}\] 
are the zeroes of $U_{\ell-1}(\tfrac{1}{2x})$. It is known that the zeroes of Chebyshev polynomials of the second kind are when $1/2x \in \{ \cos(\tfrac{n\pi}{\ell})$ : $n = 1,...,\ell-1\}$. So, the roots are $x$ such that $|x| = 1/|2\cos(\tfrac{n\pi}{k+1})| > 1/2$ for all $n$, as follows directly from the expression.

Combining these facts, we see that the growth of $b_n$ is determined by the behavior of the first term near $x=1/2$, and the second term at $x=-1/2$. A direct calculation using $U_k(1) = k+1$ shows that at $x=1/2$, the term in front of $\sum_{n\geq 0} a_n x^n $ is $ (\ell+1)/2\ell$. At $x=-1/2$, the front of the second term gives us $1/2\ell$. 

Then, the growth of our coefficients is given via 
\[b_n \sim \big(\tfrac{\ell+1+(-1)^{n+1}}{2\ell}\big) a_n.\]
Moreover, 
$|b_{n}-b_{n}^{a}|/(n^{-1/2}\cdot 2^n)\in O(\tfrac1n)$ follows identically as for $a_n$, by adapting the proof of \autoref{T:AlmostMain}, so the claim follows.
\end{proof}

\section{Background}\label{S:Background}

We assume that the reader has some background on quantum groups and their representations (in particular, tilting modules). We refer to \cite{AnPoWe-representation-qalgebras,Do-q-schur,AnStTu-cellular-tilting} for some details. Our ground field throughout is $\C$, and we regard the representation theory of the standard Lie groups and Lie algebras as well understood, in the sense that our goal is to reduce other representation theories to it.

Fix an element $q\in\C\setminus\{0\}$ and let $[k]=(q^{k}-q^{-k})/(q-q^{-1})$ for $k\in\Z_{>0}$ be the usual quantum number. Let $\ell\in\Z_{>0}\cup\{\infty\}$ minimal with respect to $[\ell]=0$. The case $\ell=\infty$ is the generic case or 
characteristic zero case.

\begin{Example}
If $q=\pm \exp(2\pi i/4)$ (the imaginary unit and its conjugate), then $\ell=2$ since $[1]=1,[2]=0$. For $q=\pm \exp(2\pi i/3)$ we have $\ell=3$ since $[1]=1,[2]=-1,[3]=0$. In general, $\ell=\infty$ unless $q$ is a root of unity, and if $q$ is a root of unity of order $k$, then $\ell=k$ if $k$ is odd, and $\ell=k/2$ if $k$ is even.
\end{Example}

For a semisimple complex Lie algebra $\mathfrak{g}=\mathfrak{g}(\C)$, let $U_q=U_q(\mathfrak{g})$ be the divided power quantum group at $q$, and work in the category $\mathbf{Rep}(U_q)$ of finite dimensional $U_q(\mathfrak{g})$-modules (of type 1). This category is a $\C$-linear abelian braided pivotal category.

\begin{Remark}
A crucial difference between working with $\C$ and fields of prime characteristic is that $\mathbf{Rep}(U_q)$ does not have enough projectives in prime characteristic.
\end{Remark}

Let $X^+$ denote dominant integral $\mathfrak{g}$-weights. Then $\mathbf{Rep}(U_q)$ has the following important indecomposable objects, indexed by $X^+$ (serving as the highest weight):
\begin{enumerate}

\item The simple module $L(\lambda)$ and its (indecomposable) projective cover $P(\lambda)$.

\item The indecomposable tilting module $T(\lambda)$.

\item The Weyl module $\Delta(\lambda)$ and its dual 
$\nabla(\lambda)$.

\end{enumerate}
For completeness, projective = injective in $\mathbf{Rep}(U_q)$, and $P(\lambda)$ is also the injective hull of $L(\lambda)$.

If $\ell=\infty$, we have $\Delta(\lambda)\cong\nabla(\lambda)\cong L(\lambda)\cong T(\lambda)\cong P(\lambda)$, and they have the same character as for $\mathfrak{g}$. In fact, 
$\mathbf{Rep}(U_q)\cong\mathbf{Rep}(\mathfrak{g})$ as $\C$-linear abelian categories (the latter is the category of finite dimensional $\mathfrak{g}$-modules). Although these categories are not equivalent as monoidal categories, their Grothendieck rings are isomorphic, i.e. the objects have the same combinatorial decomposition and multiplication structures. For $\mathfrak{g}=\mathfrak{sl}_N$ we can also alternatively work with the group $SLN=\mathrm{SL}_N(\C)$ instead of the Lie algebra. We will not distinguish between the generic quantum group and the Lie algebra or group, so its representation theory counts as well-known.

\begin{Example}\label{E:DigitGame}
Since we discuss SL2 in some details, cf. \autoref{S:SL2}, \autoref{S:SL2one} and \autoref{S:SL2two}, let us be very explicit in this case. See e.g. \cite{AnTu-tilting,SuTuWeZh-mixed-tilting} for details.

Here we have $X^+=\N$, so all above modules are indexed by natural numbers. In this case, the $(k+1)$ dimensional $\nabla(k)$, when considered for the group SL2, is the $k$th symmetric power of the defining two dimensional representation $\nabla(1)$, and $\Delta(k)$ is its dual (of the same dimension and character).

Write $[a,b]_\ell=a\cdot\ell+b$ for $a\in\N,b\in\{0,\dots,\ell-1\}$.
Let us now distinguish three cases:
\begin{enumerate}
\item The fundamental alcove is the case $a=0$.
\item The projective cone is the case $a\neq 0$.
\item The wall is the case $a\neq0,b=0$. 
\end{enumerate}
Let $X(k)$ be any of the above modules. The crucial number for all of them is $k+1$ (``off-by-one'') which we express as $[a,b]_\ell$.

Then we have:
\begin{enumerate}[label=$\bullet$]
\item $T(k)$ is simple $\Leftrightarrow$ $T(k)$ is a Weyl module $\Leftrightarrow$ $T(k)$ is a dual Weyl module $\Leftrightarrow$ $k+1$ is in the fundamental alcove or on the wall.
\item $T(k)$ is projective $\Leftrightarrow$ $k+1$ is in the projective cone.
\end{enumerate}

The $\Delta$ and equivalently $\nabla$ factors of $T(k)$ can be computed as follows. If $k+1=[a,b]_\ell=a\cdot\ell+b$ for $a\neq 0,b\in\{1,\dots,\ell-1\}$, then $T(k)$ has factors 
$\Delta(k)$ and $\Delta(j)$ where $j+1=[a,-b]_\ell=a\cdot\ell-b$. If $a=0$ or $b=0$, then $T(k)\cong\Delta(k)$. This is the flipping digits game from \cite{TuWe-quiver-tilting,SuTuWeZh-mixed-tilting}, and it implies that $\dim_\C T(k)=a\cdot\ell=k+1$ if $a\neq 0,b=0$ (on the wall), and $\dim_\C T(k)=2a\cdot\ell$ if $a\neq 0,b\neq 0$ (projective cone), and $\dim_\C T(k)=k+1$ for $a=0$ (fundamental alcove).
\end{Example}

We now give a cheat-sheet for working with tilting modules (all standard facts):
\begin{enumerate}

\item $\Delta(\lambda)$ appears $m$ times in a filtration of a tilting module $T$ if and only if 
$\nabla(\lambda)$ appears $m$. So one can use them essentially interchangeably.

\item In general, the modules $L(\lambda)$ do not have characteristic zero (non-quantum) counterparts. However, $\Delta(\lambda)$, $\nabla(\lambda)$ and $T(\lambda)$ do:

\begin{enumerate}

\item For $\Delta(\lambda)$ and $\nabla(\lambda)$ there exist a unique simple $\mathfrak{g}$-representation $\Delta_\C(\lambda)=\nabla_\C(\lambda)$ of the same character.

\item Assume that $T(\lambda)$ has $\Delta$-factors $\{\Delta(\lambda_i)\}$. Then $T_\C(\lambda)=\bigoplus_i\Delta_\C(\lambda_i)$ is a $\mathfrak{g}$-representation of the same character.

\end{enumerate}

\item By taking direct sums, the previous point gives a characteristic zero counterpart $T_\C$ for any tilting module $T$.
In particular, for every tilting module $T$, the number of direct summands in $T$ is less than (or equal to) the number of direct summands of $T_\C$.

\item A tilting module is uniquely determined by its $\Delta$-factors.

\end{enumerate}

\begin{Example}
With the notation of \autoref{E:DigitGame}, for $\ell=3$ the tilting module $T(3=4-1)$ has $\Delta$-factors $\Delta(3=4-1)$ and $\Delta(1=2-1)$ since $4=[1,1]_\ell$ and $2=[1,-1]_\ell$. So $T(3)^{\otimes 2}$ has $\Delta$-factors $\Delta(3)^{\otimes 2}\cong \Delta(6)\oplus\Delta(4)\oplus\Delta(2)\oplus\Delta(0)$, twice $\Delta(3)\otimes\Delta(1)\cong\Delta(4)\oplus\Delta(2)$, and $\Delta(1)^{\otimes 2}\cong\Delta(2)\oplus\Delta(0)$. Thus, by collecting $\Delta$-factors, we get $T_\C(3)^{\otimes 2}\cong \Delta(6)\oplus \Delta(4)^{\oplus 3}\oplus \Delta(2)^{\oplus 4}\oplus\Delta(0)^{\oplus 2}$ and
$T(3)^{\otimes 2}\cong T(6)\oplus T(4)^{\oplus 3}\oplus T(2)\oplus T(0)^{\oplus 2}$.
\end{Example}

We will summarize some alcove
combinatorics, which will depend on $\ell$. Let $I$ be the indexing set for the vertices of the Dynkin diagram associated to $\mathfrak{g}$. We use the usual notation for the underlying (positive) Euclidean space $E$ (or $E^+$), the roots $\alpha_i$, coroots $\alpha_i^{\vee}$, fundamental weights $\omega_i$ etc. Let $\rho$ be the half-sum of the positive roots.
We also use the usual pairings, as in \autoref{E:Combinatorics}.

\begin{Example}\label{E:Combinatorics}
Explicitly, for SLN we can identify everything in
the above statistic so that: $I=\{1,\dots,N-1\}$,
$E=\R^{N-1}$ with standard basis $\{e_i|i\in I\}$ and pairing $(e_i,e_j)=\delta_{ij}$,$E^+=\R_{\geq 0}^{N-1}$, $\alpha_i=e_i-e_{i+1}$, $\alpha_i^{\vee}$ determined by $\langle\alpha_i,\alpha_j^{\vee}\rangle=(\alpha_i,\alpha_j)=i,j$th entry of the Cartan matrix. The fundamental weights are $\omega_i=\alpha_1+\dots+\alpha_i$, and we can and will write $\lambda\in X^+\cong\N^{N-1}$ uniquely in terms of fundamental weights as $(\lambda_1,\dots,\lambda_{N-1})\in\N^{N-1}$. In this notation, $\rho=(1,\dots,1)$.
\end{Example}

The affine Weyl group $W_{\ell}$ is
the group generated by reflections $\{s_{i,r}\mid i\in I,r\in\Z\}$ in the
affine hyperplanes $H^{\ell}_{i,r}=\{x\in E\mid\langle
x+\rho,\alpha_i^{\vee}\rangle=\ell r\}$. The
group $W_{\ell}$ acts on $E$ via the ($\ell$-scaled) dot-action which in coordinates is
$s_{i,r}\bullet x=x+(r\ell-\langle x,\alpha_i^{\vee}\rangle-2)\alpha_i$, the orbits
under this action intersected with $X^+$ are (sometimes) called blocks.

The
fundamental alcove is given by $A_0=\{x\in E^+\mid 0<\langle
x+\rho,\alpha_i^{\vee}\rangle<\ell, i\in I\}$, and in general an
alcove $A\subset E^+$ is a connected component of
$\big(E\setminus\bigcup_{i,r}H^{\ell}_{i,r}-\rho\big)\cap E^+$. It is easy to
see that such alcoves consist of elements in $E^+$ satisfying strict
inequalities similarly to the fundamental alcove, and the upper closure
$\widehat{A}$ of such an alcove is then defined by replacing the second $<$ by
$\leq$.

Moreover, a $\lambda\in X^+$ is called $k$ singular if its dot-action stabilizer is of the corresponding size $k$. (The $0$ singular weights are often called regular.) The property of being $k$ singular is shared by all dominant weights in the same block.
Further, in each upper closure of an alcove there is exactly one $\lambda\in X^+$
per block, and if we have a particular block in mind, we identify this element $\lambda\in X^+$ with its alcove.

Finally,
$X^+$ is partitioned into three sets $I_0=\widehat{A}_0\cap X^+$,
$I_p=(\ell-1)\rho+\N\sum_{i\in I}\omega_i$ and
$I_b=X^+\setminus(I_0\cup I_p)$, the latter of which we
call the boundary case; the set $I_0$ contains the point
$0\in E$, while $I_p$ contains the so-called Steinberg point
$(\ell-1)\rho$.

\begin{figure}[ht]
\begin{tikzpicture}[anchorbase]
\node at (0,0) {\includegraphics[height=10cm]{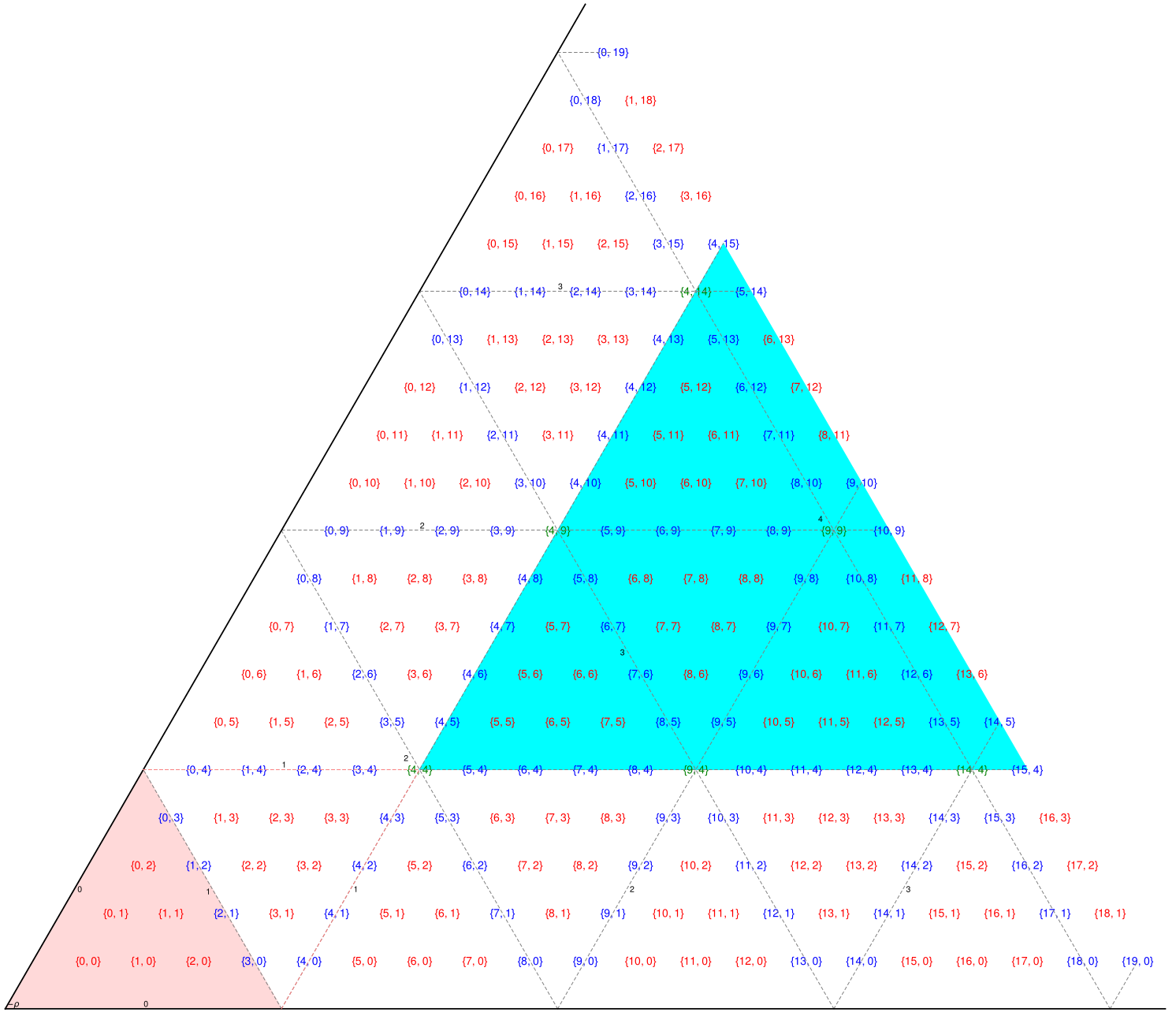}};
\end{tikzpicture}. 
\caption{The
dominant $SL3$ weights $(m,n)\in X^+$, distinguished by color depending on their singularity, arranged in $\ell=5$ scaled alcoves. The
alcove containing $I_0$ and the set $\partial\widehat{A}_1$ are illustrated in pink, those containing $I_p$ are illustrated in cyan, and the ones containing $I_b$ form the remaining white collar.}
\label{figure:sl3-combinatorics}
\end{figure}

\autoref{figure:sl3-combinatorics} summarizes these conventions for SL3.

\begin{Example}
The case of SL2 is special, since $I_b=\emptyset$. Moreover, $I_0=\{0,\dots,\ell-2\}$ and $I_p=\Z_{\geq(\ell-1)}$. The Steinberg point is $\ell-1$.
\end{Example}

For $\ell\geq h$ (where $h$ is the Coxeter number associated to $\mathfrak{g}$), the simple $U_q$-tilting modules are the ones 
with defining weight in $I_0$ or with maximally singular defining weight, and the projectives are the ones with defining weight in $I_p$.
An important simple $U_q$-tilting module is the so-called Steinberg module $\mathrm{ST}=T\big((\ell-1)\rho\big)$ of dimension $\dim_\C\mathrm{ST}=\ell^{\#R^+}$, where $\#R^+$ is the number of positive roots associated to $\mathfrak{g}$.

\section{SL2 and quantum SL2}\label{S:SL2one}

We start with the harvest of \autoref{S:SL2}. Let $U_q=U_q(\mathfrak{sl}_2)$ defined as in \autoref{S:Background}.

\begin{Theorem}\label{T:SL2oneMainTwo}
For $b_n=b_n\big(T(1)\big)$, we have
\begin{gather*}
b_n\sim
\begin{cases}
\tfrac{\ell+1}{2\ell} \cdot a_n^a
& \ell \text{ odd},\\[2ex]
\tfrac{\ell+1+(-1)^{n+1}}{2\ell}\cdot a_n^a, & \ell \text{ even, }
\end{cases}
,
\end{gather*}
with variance in $O(n^{-3/2}\cdot 2^n)$.
Moreover, for any tilting $U_q$-representation $V=T(k-1)$ and $b_n=b_n\big(T(k-1)\big)$ we have
\begin{gather*}
\tfrac{1}{2}\cdot\sqrt{\tfrac{6}{(k^2-1)\pi}}\cdot n^{-1/2}\cdot(\dim_\C V)^n\leq b_n\leq \sqrt{\tfrac{6}{(k^2-1)\pi}}\cdot n^{-1/2}\cdot(\dim_\C V)^n,\quad\text{for $n\gg 0$}.
\end{gather*}
\end{Theorem}

\begin{proof}
We comment now on the numbers defined in \autoref{Eq:Recursion}.

Let $\C^2\cong\nabla_\C(1)$ be the defining (classical) representation of $SL2=\mathrm{SL}_2(\C)$, where a matrix in SL2 acts as itself on $\nabla(1)$. Extending the action diagonally, we get an associated action on $\nabla_\C(1)^{\otimes n}$. Then $a_n$ is the number of indecomposable summands of $\nabla_\C(1)^{\otimes n}$ (equivalently, the number of simple summands). The numbers $a_{n,m}$ count the number of individual summands $\mathrm{Sym}^m\C^2\cong\nabla_\C(m)$ (these are precisely the indecomposables of SL2) in $\nabla_\C(1)^{\otimes n}$.

The numbers $b_n$ are a quantum version of this. In more detail, 
let $T(1)$ be the two dimensional indecomposable tilting module over $U_q$ that dequantizes (in a precise sense) to $\nabla_\C(1)$. Let $T(m)$ be the tilting modules with the same highest weight as $\mathrm{Sym}^m\C^2$. Then $b_{n,m}$ is the number of appearances of $T(m)$ in $V^{\otimes n}$, and $b_n$ = $\sum_{m \leq n} b_{n,m}$ is the total number of indecomposable summands of $V^{\otimes n}$. The recursion of $b_{n,m}$ in \autoref{Eq:Recursion} is \cite[Corollary 3.5]{An-simple-tl}.

The claims for $T(1)$ then follow directly from \autoref{T:MainSL2}. 

Now assume that $V$ is tilting with highest weight $k-1$. Let $V_\C$ be the characteristic zero counterpart of it, namely the 
$\mathfrak{sl}_2$-representation that is the direct sum of the Weyl modules of $V$. By \cite{Bi-asymptotic-lie,CoEtOs-growth} we have
\begin{gather*}
b_n(V_\C)\sim \sqrt{\tfrac{6}{(k^2-1)\pi}}\cdot n^{-1/2}\cdot(\dim_\C V)^n.
\end{gather*}
Not that $b_n(V_\C)$ counts the Weyl factors of $V^{\otimes n}$ since the Weyl modules are direct summands for $\mathfrak{sl}_2$ over $\C$.

\begin{Lemma}
Every indecomposable tilting module for $U_q$ has one or two Weyl factors.
\end{Lemma}

\begin{proof}
Directly from \autoref{E:DigitGame}.
\end{proof}

Now, this shows
\begin{gather*}
2^{-1}\cdot b_n(V_\C)\leq b_n\leq b_n(V_\C),
\end{gather*}
and the claim follows by taking everything together.
\end{proof}

\begin{Remark}
The results of \cite{LaPoReSo-growth-qgroup} are somewhat orthogonal: They compute exact multiplicities in $T(1)^{\otimes n}$ at even roots of unity, via a lattice path model, refining parts of \autoref{T:SL2oneMainTwo}. Conversely, our contribution is broader in scope: it treats all roots of unity, all tilting $U_q$-modules and extends to general type (see \autoref{S:General}), where the dominant behavior depends only on the dimension and the root system rather than the specific module. In short, \cite{LaPoReSo-growth-qgroup} gives sharp SL2 multiplicities using a lattice path model; we provide a uniform, representation-agnostic picture using random walks and tilting-module methods.
\end{Remark}

Back to $V=T(1)$, we add the following representation-theoretical analog of \autoref{T:AlmostMain}. For $n\in\N$, let $w_n$ be the number of indecomposable summands of $V^{\otimes n}$ whose highest weight lies on the wall, i.e.\ the number of indecomposable summands $T(m)$ with $m+1\equiv 0\bmod\ell$.

\begin{Theorem}\label{T:WallSummands}
We have 
\[
w_n\sim
w_n^a:=
\begin{cases}
\tfrac{1}{\ell}\cdot a_n^a, & \ell \text{ odd},\\
0, & \ell \text{ even and }n\text{ even},\\
\tfrac{2}{\ell}\cdot a_n^a, & \ell \text{ even and }n\text{ odd},
\end{cases}
\quad
\lim_{n\to\infty}\tfrac{w_n}{b_n}=
\begin{cases}
\tfrac{2}{\ell+1}, & \ell \text{ odd},\\
0, & \ell \text{ even, along even }n,\\
\tfrac{4}{\ell+2}, & \ell \text{ even, along odd }n,
\end{cases}
\]
and moreover $|w_n-w_n^a|/(n^{-1/2}\cdot 2^n)\in O(\tfrac1n)$
\end{Theorem}

\begin{proof}
The number $w_n$ of wall summands in $T(1)^{\otimes n}$ is
$w_n=c_n^{(\ell-1)}$. Applying \autoref{T:AlmostMain}.(b) with $r=\ell-1$ gives the stated asymptotics for $w_n$, including the $O(\tfrac1n)$ relative error term. This proves the first part, and we can directly use 
\autoref{T:MainSL2}, completing the proof.
\end{proof}

\section{General results}\label{S:General}

Let $\#R^+$ denote the number of positive roots associated with $\mathfrak{g}$, and let 
$U_q=U_q(\mathfrak{g})$ be as in \autoref{S:Background}.
The following is the main general result (generalizing \cite[Section 4]{CoOsTu-growth}) and will be proven in the remainder of this section. Using the usual Bachmann--Landau notation, let $g\in\Theta(f)$ denote ($g\in O(f)$ and $f\in O(g)$).

\begin{Theorem}\label{T:MainTwo}
For any tilting $U_q$-representation $T$ and $b_n=b_n(T)$, we have
\begin{gather*}
b_n\in\Theta\big(n^{-\#R^+/2}\cdot(\dim_\C T)^n\big).
\end{gather*}
\end{Theorem}

\begin{Example}
For SL2 and its defining representation, \autoref{T:MainTwo} gives
\begin{gather*}
b_n\in\Theta\big(n^{-1/2}\cdot 2^n\big),
\end{gather*}
which is weaker than \autoref{T:MainSL2} (and \autoref{T:SL2oneMainTwo}).
\end{Example}

We start with the following result.

\begin{Lemma}\label{L:ProjTiltingBasics}
In $\mathbf{Rep}(U_q)$ we have:
\begin{enumerate}
\item Every indecomposable projective object of $\mathbf{Rep}(U_q)$ is a tilting module.
\item The full subcategory of projective objects is a two-sided tensor ideal: if $P$ is projective and $M$ is any object, then both $P\otimes M$ and $M\otimes P$ are projective.
\item The class of tilting modules is closed under tensor products: if $T,U$ are tilting, then $T\otimes U$ is tilting.
\item In (c), $T\otimes U$ has the same $\Delta$-factors as $T_\C\otimes U_\C$.
\end{enumerate}
\end{Lemma}

\begin{proof}
Standard arguments: (a) Steinberg generates projectives \cite[Section 9]{AnPoWe-representation-qalgebras}; (b) rigidity implies tensor-ideal property \cite[Proposition 4.2.12]{EtGeNiOs-tensor-categories}; (c) tensor product of modules with Weyl/dual-Weyl filtrations again has both filtration \cite{Pa-tilting-tensor}, and then (d) follows from the construction in \cite{Pa-tilting-tensor}.
\end{proof}

\begin{Lemma}\label{L:WeylBound}
The number of $\Delta$-factors in a tensor product 
$\Delta(\mu)\otimes\Delta(\nu)$ of Weyl modules is at most $\min\{\dim_\C \Delta(\nu),\dim_\C \Delta(\mu)\}$.
\end{Lemma}

\begin{proof}
By viewing the bimodule $\Delta(\mu)\otimes\Delta(\nu)$ as a left or right module, where the other side of the tensor product is a multiplicity space.
\end{proof}

Let $h$ denote the Coxeter number associated to $\mathfrak{g}$ (i.e. the order of any Coxeter element in the associated Weyl group). For this well-known (and easy to find) statistic associated to $\mathfrak{g}$ we get:

\begin{Lemma}\label{L:KL}
If $\ell\geq h$, then $\big(T(\lambda),\Delta(\mu)\big)$ (the Weyl multiplicities) are governed by (parabolic) antispherical Kazhdan--Lusztig polynomials.
\end{Lemma}

\begin{proof}
This is a main result of \cite{So-tilting-a,So-tilting-b}. (These papers only treat the case of odd $\ell$, but with \cite{An-strong-link}, the results also follow for even $\ell$.)
\end{proof}

Let $(P(\lambda),\Delta)=\sum_{\mu\in X^+}\big(P(\lambda),\Delta(\mu)\big)=\sum_{\mu\in X^+}\big(T(\lambda^{\prime}),\Delta(\mu)\big)$ be the number of Weyl factors in $P(\lambda)$. (All but finitely many summands are zero, so $(P(\lambda),\Delta)\in\N$.)

\begin{Lemma}\label{L:Bound}
We have
\[
(P(\lambda),\Delta)\leq\ell^{\,\#R^+}.
\]
Moreover, if $\ell\geq h$, then we have
\[
(P(\lambda),\Delta)\leq h^{\,\#R^+}.
\]
\end{Lemma}

\begin{proof}
Write $\lambda=\ell\lambda^1+\lambda^0$ with $\lambda^0$ restricted.
By Steinberg tensor product one has
\[
P(\lambda) \;\cong\; P(\lambda^0)\otimes \Delta(\lambda^1)^{[q]},
\]
where the exponent denotes the Frobenius--Lusztig twist.
In the restricted case $\lambda=\lambda^0$, the module $P(\lambda^0)$ is a 
direct summand of
\[
\mathrm{ST} \otimes \Delta_q\!\big((\ell-1)\rho+w_0\!\bullet \lambda^0\big),
\]
where $\mathrm{ST}=T\big((\ell-1)\rho\big)$ is the Steinberg module (with defining weight the Steinberg point) as above, which is simple and projective-injective of dimension $\dim_\C\mathrm{ST}=\ell^{\#R^+}$. This follows from the result in and leading to \cite[\S 8–9]{AnPoWe-representation-qalgebras}.

Applying \autoref{L:WeylBound} to the summand above shows that
\[
\big(P(\lambda^0),\Delta\big) \leq \dim_\C\mathrm{ST}=\ell^{\#R^+}.
\]
The same argument applies in general $\lambda=\ell\lambda^1+\lambda^0$, 
since then $P(\lambda)$ is a summand of
\[
\Delta\big((\ell-1)\rho+\ell\lambda^1\big)\otimes 
\Delta\big((\ell-1)\rho+w_0\!\bullet\lambda^0\big),
\]
and $\dim_\C \Delta((\ell-1)\rho+w_0\!\bullet\lambda^0)\leq \dim_\C \mathrm{ST}$ which follows since $\mathrm{ST}$ is a Weyl module and hence, the classical Weyl dimension formula applies on both sides of the inequality.
Hence the uniform bound follows.

Finally, for $\ell\geq h$, this is essentially independent of $\ell$ by \autoref{L:KL} (the multiplicities are governed by antispherical Kazhdan--Lusztig polynomials which stabilize with $\ell$), so we can take the smallest possible $\ell$, which is $h$.
\end{proof}

Note that the above implies that every projective $U_q$-module
$P$ has an associated characteristic zero (for $\mathfrak{g}$) analog $P_\C$ with the same character.

\begin{Lemma}\label{L1}
For any tilting $U_q(\mathfrak{g})$-module $T$ and any projective $P$, the number $b(P\otimes T)$ of indecomposable summands of $P\otimes T$ can be bounded:
\begin{gather*}
\text{general: }
\ell^{-\#R^+}\cdot b(P_\C\otimes T_\C)\leq b(P\otimes T)\leq b(P_\C\otimes T_\C).,
\\
\text{for $\ell\geq h$: }
h^{-\#R^+}\cdot b(P_\C\otimes T_\C)\leq b(P\otimes T)\leq b(P_\C\otimes T_\C).
\end{gather*}
\end{Lemma}

\begin{proof}
The upper bound comes from comparing filtrations with characteristic zero.  
The lower bound comes from the uniform factor $\ell^{\#R^+}$ or $h^{\#R^+}$ in \autoref{L:Bound}.
\end{proof}

\begin{Lemma}\label{L:FBC}
Let $T$ be a $U_q$-tilting module that contains at least one nontrivial summand. Then, for any $\lambda\in X^+$, there exists $n\in\Z$ such that $P(\lambda)$ is a direct summand of $T^{\otimes n}$. Moreover, the projective cell is a FBC for the fusion graph associated to $T$, and it is the only basic class and the only final class.
\end{Lemma}

\begin{proof}
The first part, and that the projective cell is a FBC and the only final class, follows from the proof of \cite[Proposition 4.23]{LaTuVa-growth-pfdim-inf} and \autoref{L:ProjTiltingBasics}.

To show that there are no other basic classes, let $\Gamma^\prime(T)$ be the fusion graph associated with $T$ where we identify strongly connected components to vertices that are connected by an edge if tensoring with $T$ gets us from one strongly connected component to the other. It follows from \cite{Os-tensor-ideals-tilting} that the graph $\Gamma^\prime(T)$ is finite, and the projective cell is a FBC in $\Gamma^\prime(T)$. Then \cite[Theorem 8]{LaTuVa-growth-pfdim} gives that there cannot be another basic class (since otherwise the growth would be strictly larger than in characteristic zero).
\end{proof}

\begin{Lemma}\label{L2}
Let $T$ be a $U_q$-tilting module that contains at least one nontrivial summand. Then the associated growth problem $b_n=b(T^{\otimes n})$ is determined on the ideal of projectives in $\mathbf{Rep}(U_q)$.
\end{Lemma}

\begin{proof}
Directly from \autoref{L:FBC}.
\end{proof}

\begin{proof}(Proof of \autoref{T:MainTwo}.)
Combine the characteristic zero growth of $b(T_\C^{\otimes n})$ (which is known by \cite{Bi-asymptotic-lie}) with \autoref{L1} and \autoref{L2}.
\end{proof}

\section{Better bounds in rank 2}\label{S:Better}

Let us assume that $\ell>2$ for SL3 (type A2), 
$\ell>4$ for SO5 and SP4 (types B2 and C2), and 
$\ell>6$ for G2. The following generalizes \cite[Proposition 4.2]{CoOsTu-growth}.

\begin{Theorem}
We can improve the lower bound $a^{-1}=h^{-\#R^+}$ in \autoref{L1} as follows:
\begin{enumerate}
\item For type A2, the lower bound can be chosen to be $a=12$ (instead of $3^3=27$).
\item For types B2 and C2, the lower bound can be chosen to be $a=32$ (instead of $4^4=256$).
\item For type G2, the lower bound can be chosen to be $a=348$ (instead of $6^6=46656$).
\end{enumerate}
\end{Theorem}

\begin{proof}
We begin with a well-known lemma, that, as far as we are aware, is not in the literature 
in the formulation adapted to our needs. So we give a short proof.

Fix a regular dominant weight $\lambda$ and a singular dominant weight $\lambda'$ (of type $J$), and let
\[
\Theta=\Theta_{\lambda\to\lambda'}
\]
be translation onto the $J$-wall.

\begin{Lemma}\label{L:Translation}
Then for all tilting $U_q$-modules $T$,
\[
\ell_\Delta(\Theta T)\ \le\ \ell_\Delta(T),
\]
with strict inequality whenever some $\Delta_q(\nu)$ occurring in a $\Delta$-filtration of $M$ lies outside the $\lambda'$-block. In particular, for the indecomposable tilting $T(\lambda)$ and any indecomposable summand $U$ of $\Theta T(\lambda)$ one has
\[
\ell_\Delta(U)\ \le\ \ell_\Delta\big(T(\lambda)\big),
\]
and $\ell_\Delta(\Theta T(\lambda))<\ell_\Delta(T(\lambda))$ provided $J\neq\varnothing$.
\end{Lemma}

\begin{proof}
Translation onto a wall is exact, preserves tiltings and $\Delta$–filtrations, and is given by tensoring with a finite dimensional tilting followed by projection to the $\lambda'$-block. If
$0=M_0\subset\cdots\subset M_r=M$ with $M_i/M_{i-1}\cong \Delta_q(\mu_i)$, then exactness yields a filtration of $\Theta M$ whose successive quotients are $\Theta\Delta_q(\mu_i)$, each either $0$ or a Weyl module in the $\lambda'$-block. Thus no new $\Delta$–factors appear and some may vanish, giving $\ell_\Delta(\Theta M)\le r=\ell_\Delta(M)$, with strictness exactly when one of the $\Theta\Delta_q(\mu_i)$ is $0$. Additivity of $\ell_\Delta$ over direct sums gives the statement for summands.
\end{proof}

With \autoref{L:Translation} in place, we simply observe that there are only finitely many antispherical KL polynomials for the indecomposable projective tilting modules, and the largest ones are of the respective sizes \cite{St-diplom}. Now, we use the fact that projective titling modules form the unique FBC in $\mathbf{Rep}(U_q)$ via \autoref{L2}.
\end{proof}

\section{SL2 in the mixed case}\label{S:SL2two}

We now consider SL2 with its defining representation again, but for a field $\K$ of arbitrary characteristic $p$ (with characteristic zero considered as $p=\infty$).
The below does not depend on the field, but only on the pair $(p,\ell)$.

Let $b_n^{\ell=p}$ denote the case of $\ell=p$.
The case $b_n^{\ell=p}$ was discussed in \cite{CoEtOsTu-growth-fractal}, and the following partially generalizes this (\cite{CoEtOsTu-growth-fractal} provides more information about the special case $\ell=p$).

\begin{Theorem}\label{T:Mixed}
We have
\begin{gather*}
b_n\in\Theta(b_n^{\ell=p}).
\end{gather*}
\end{Theorem}

\begin{proof}
The case $p=\infty$ was discussed in \autoref{S:SL2}, and for $\ell=\infty$ the situation is semisimple (so it is the same as for $\mathrm{SL}_2(\C)$), so assume $p,\ell\neq\infty$.

Following \cite{CoEtOsTu-growth-fractal}, for $w\in\C$ with $|w|<1$ let
\begin{gather*}
F(w)=\sum_{n\geq 0}b_n(w+w^{-1})^{-n-1},
\quad
F(w)^{\ell=p}=\sum_{n\geq 0}b_n^{\ell=p}(w+w^{-1})^{-n-1}
\end{gather*}
Repeating the arguments in \cite{CoEtOsTu-growth-fractal}, with the flipping digits game as explained in 
\autoref{E:DigitGame}, one can see that
\begin{gather*}
F(w)=\tfrac{w(1-w^{\ell-1})(1-w)^{-1}(1+w^{\ell})^{-1}}{w(1-w^{p-1})(1-w)^{-1}(1+w^{p})^{-1}}\cdot F(w)^{\ell=p}
=
\tfrac{(1-w^{\ell-1})(1+w^{p})}{(1-w^{p-1})(1+w^{\ell})}\cdot F(w)^{\ell=p}
.
\end{gather*}
For the limit $w\uparrow1$, the factor in front is $\tfrac{\ell-1}{p-1}$, so is not important for a 
capital Theta analysis. The result follows.
\end{proof}

For SL2, note that \autoref{T:Mixed} generalizes \autoref{T:MainTwo} and expresses the quantum case as the $q=1$ case over the same field. One would expect something similar for general quantum groups.


\end{document}